

%
\magnification1200
\pretolerance=100
\tolerance=200
\hbadness=1000
\vbadness=1000
\linepenalty=10
\hyphenpenalty=50
\exhyphenpenalty=50
\binoppenalty=700
\relpenalty=500
\clubpenalty=5000
\widowpenalty=5000
\displaywidowpenalty=50
\brokenpenalty=100
\predisplaypenalty=7000
\postdisplaypenalty=0
\interlinepenalty=10
\doublehyphendemerits=10000
\finalhyphendemerits=10000
\adjdemerits=160000
\uchyph=1
\delimiterfactor=901
\hfuzz=0.1pt
\vfuzz=0.1pt
\overfullrule=5pt
\hsize=146 true mm
\vsize=8.9 true in
\maxdepth=4pt
\delimitershortfall=.5pt
\nulldelimiterspace=1.2pt
\scriptspace=.5pt
\normallineskiplimit=.5pt
\mathsurround=0pt
\parindent=20pt
\catcode`\_=11
\catcode`\_=8
\normalbaselineskip=12pt
\normallineskip=1pt plus .5 pt minus .5 pt
\parskip=6pt plus 3pt minus 3pt
\abovedisplayskip = 12pt plus 5pt minus 5pt
\abovedisplayshortskip = 1pt plus 4pt
\belowdisplayskip = 12pt plus 5pt minus 5pt
\belowdisplayshortskip = 7pt plus 5pt
\normalbaselines
\smallskipamount=\parskip
 \medskipamount=2\parskip
 \bigskipamount=3\parskip
\jot=3pt
%
%
\def\ref#1{\par\noindent\hangindent2\parindent
 \hbox to 2\parindent{#1\hfil}\ignorespaces}
%
%
\font\tenss=cmss10
\font\sevenss=cmss8 at 7pt
\font\fivess=cmss8 at 5pt
\newfam\ssfam %
\textfont\ssfam=\tenss
\scriptfont\ssfam=\sevenss
\scriptscriptfont\ssfam=\fivess
%
%
%
%
%
%
%
%
%
\catcode`\_=11
\def\suf_fix{}
\def\scaled_rm_box#1{%
 \relax
 \ifmmode
   \mathchoice
    {\hbox{\tenrm #1}}%
    {\hbox{\tenrm #1}}%
    {\hbox{\sevenrm #1}}%
    {\hbox{\fiverm #1}}%
 \else
  \hbox{\tenrm #1}%
 \fi}
\def\suf_fix_def#1#2{\expandafter\def\csname#1\suf_fix\endcsname{#2}}
\def\I_Buchstabe#1#2#3{%
 \suf_fix_def{#1}{\scaled_rm_box{I\hskip-0.#2#3em #1}}
}
\def\rule_Buchstabe#1#2#3#4{%
 \suf_fix_def{#1}{%
  \scaled_rm_box{%
   \hbox{%
    #1%
    \hskip-0.#2em%
    \lower-0.#3ex\hbox{\vrule height1.#4ex width0.07em }%
   }%
   \hskip0.50em%
  }%
 }%
}
\I_Buchstabe B22
\rule_Buchstabe C51{34}
\I_Buchstabe D22
\I_Buchstabe E22
\I_Buchstabe F22
\rule_Buchstabe G{525}{081}4
\I_Buchstabe H22
\I_Buchstabe I20
\I_Buchstabe K22
\I_Buchstabe L20
\I_Buchstabe M{20em }{I\hskip-0.35}
\I_Buchstabe N{20em }{I\hskip-0.35}
\rule_Buchstabe O{525}{095}{45}
\I_Buchstabe P20
\rule_Buchstabe Q{525}{097}{47}
\I_Buchstabe R21 
\rule_Buchstabe U{45}{02}{54}
\suf_fix_def{Z}{\scaled_rm_box{Z\hskip-0.38em Z}}
\catcode`\"=12
\newcount\math_char_code
\def\suf_fix_math_chars_def#1{%
 \ifcat#1A
  \expandafter\math_char_code\expandafter=\suf_fix_fam
  \multiply\math_char_code by 256
  \advance\math_char_code by `#1
  \expandafter\mathchardef\csname#1\suf_fix\endcsname=\math_char_code
  \let\next=\suf_fix_math_chars_def
 \else
  \let\next=\relax
 \fi
 \next}
%
%
%
%
\def\font_fam_suf_fix#1#2 #3 {%
 \def\suf_fix{#2}
 \def\suf_fix_fam{#1}
 \suf_fix_math_chars_def #3.
}
\font_fam_suf_fix
 0rm
 ABCDEFGHIJKLMNOPQRSTUVWXYZabcdefghijklmnopqrstuvwxyz
\font_fam_suf_fix
 2scr
 ABCDEFGHIJKLMNOPQRSTUVWXYZ
\font_fam_suf_fix
 \slfam sl
 ABCDEFGHIJKLMNOPQRSTUVWXYZabcdefghijklmnopqrstuvwxyz
\font_fam_suf_fix
 \bffam bf
 ABCDEFGHIJKLMNOPQRSTUVWXYZabcdefghijklmnopqrstuvwxyz
\font_fam_suf_fix
 \ttfam tt
 ABCDEFGHIJKLMNOPQRSTUVWXYZabcdefghijklmnopqrstuvwxyz
\font_fam_suf_fix
 \ssfam
 ss
 ABCDEFGHIJKLMNOPQRSTUVWXYZabcdefgijklmnopqrstuwxyz
\catcode`\_=8
\def\Ndss{{\fam\ssfam I\mkern -2.5mu N}}%
\def\Qdss{{\fam\ssfam
    \mkern 3.8 mu \mathchoice%
    {\vrule height 6.5pt depth -.67pt width 1pt}%
    {\vrule height 6.5pt depth -.7pt width 1pt}%
    {\vrule height 4.55pt depth -.44pt width .7pt}%
    {\vrule height 3.25pt depth -.3pt width .5pt}%
    \mkern -5.9mu Q}}%
\def\Zdss{{\fam\ssfam Z\mkern-8.1mu Z}}%
%
%
%
%
\font\teneuf=eufm10 
\font\seveneuf=eufm7
\font\fiveeuf=eufm5
\newfam\euffam \def\euf{\fam\euffam\teneuf} 
\textfont\euffam=\teneuf \scriptfont\euffam=\seveneuf
\scriptscriptfont\euffam=\fiveeuf
\def\Afr{{\euf A}}

       \def\gfr{{\euf g}}

       \def\pfr{{\euf p}}

\parindent=0pt

\def\Hom{{\rm Hom}}
\def\Mod{{\rm Mod}}
\def\top{{\rm top}}
\def\cot{{\rm cot}}
\def\cont{{\rm cont}}
\def\comp{{\rm comp}}
\def\Ban{{\rm Ban}}
\def\coker{\mathop{\rm coker}\nolimits}
\def\ker{\mathop{\rm ker}\nolimits}
\def\im{\mathop{\rm im}\nolimits}
\def\lim{\mathop{\rm lim}}
\def\L{{\cal L}}
\def\subsetneqq{\mathop{\subset}\limits_{\neq}}
\def\dZ{{\Zdss_p}}
\def\dQ{{\Qdss_p}}
\def\Kchi{{K^{(\chi)}}}
\def\Ind{{\rm Ind}}

\centerline{\bf Banach space representations and Iwasawa theory}

\medskip

\centerline{\bf P. Schneider, J. Teitelbaum}

\bigskip

The lack of a $p$-adic Haar measure causes many methods of traditional
representation theory to break down when applied to continuous
representations of a compact $p$-adic Lie group $G$ in Banach spaces
over a given $p$-adic field $K$. For example, the abelian group
$G=\dZ$ has an enormous wealth of infinite dimensional, topologically
irreducible Banach space representations, as may be seen in the paper
by Diarra [Dia]. We therefore address the problem of finding an
additional ''finiteness'' condition on such representations that will
lead to a reasonable theory. We introduce such a condition that we
call ''admissibility''. We show that the category of all admissible
$G$-representations is reasonable -- in fact, it is abelian and of a
purely algebraic nature -- by showing that it is anti-equivalent to
the category of all finitely generated modules over a certain kind of
completed group ring $K[[G]]$.

In the first part of our paper we deal with the general
functional-analytic aspects of the problem. We first consider the
relationship between $K$-Banach spaces and compact, linearly
topologized $o$-modules where $o$ is the ring of integers in $K$. As a
special case of ideas of Schikhof [Sch], we recall that there is an
anti-equivalence between the category of $K$-Banach spaces and the
category of torsionfree, linearly compact $o$-modules, provided one
tensors the Hom-spaces in the latter category with $\Qdss$. In
addition we have to investigate how this functor relates certain
locally convex topologies on the Hom-spaces in the two categories.
This will enable us then to derive a version of this anti-equivalence
with an action of a profinite group $G$ on both sides relating
$K$-Banach space representations of $G$ and certain topological
modules for the ring $K[[G]] := K\otimes_o o[[G]]$.

Having established these topological results we assume that $G$ is a
compact $p$-adic Lie group and focus our attention on the Banach
representations of $G$ that correspond under the anti-equivalence to
finitely generated modules over the ring $K[[G]]$. We characterize
such Banach space representations intrinsically. We then show that the
theory of such "admissible" representations is purely algebraic -- one
may "forget" about topology and instead study finitely generated
modules over the noetherian ring $K[[G]]$.

As an application of our methods we determine the topological
irreducibility as well as the intertwining maps for representations of
$GL_2(\dZ)$ obtained by induction of a continuous character from the
subgroup of lower triangular matrices. Let us stress the fact that
topological irreducibility for an admissible Banach space
representation corresponds to the algebraic simplicity of the dual
$K[[G]]$-module. It is indeed the latter which we will analyze. These
results are a complement to the treatment of the locally analytic
principal series representations studied in [ST1].

\smallskip

Throughout this paper $K$ is a finite extension of $\Qdss_p$ with
ring of integers $o \subseteq K$ and absolute value $|\ |$. A
topological $o$-module is called linear-topological if the zero
element has a fundamental system of open neighbourhoods consisting
of $o$-submodules. We let
$$
\matrix{
\Mod_{\top} (o):= & \hbox{category of all Hausdorff linear-topological}\;
o\hbox{-modules}\hfill\cr & \hbox{with morphisms being all
continuous}\; o\hbox{-linear maps.}\hfill}
$$

\bigskip

{\bf 1. A duality for Banach spaces}

\smallskip

In this section we will recall a certain duality theory for
$K$-Banach spaces due to Schikhof ([Sch]). Because of the
fundamental role it will play in our later considerations and
since it is quite easy over locally compact fields we include the
proofs. We set
$$
\matrix{
\Mod^{\rm fl}_{\comp} (o):= & \hbox{the full subcategory in}\;
\Mod_\top (o)\; \hbox{of all}\hfill\cr
& \hbox{torsionfree and compact linear-topological}\hfill\cr &
o\hbox{-modules.}\hfill}
$$

\medskip

{\bf Remark 1.1:}

{\it i. An $o$-module is torsionfree if and only if it is flat;

ii. a compact linear-topological $o$-module $M$ is flat if and
only if $M\cong \prod\limits_{i\in I} o$ for some set $I$.}

Proof: i. [B-CA] Chap. I \S 2.4 Prop. 3(ii). ii. [SGA3] ${\rm
Exp.\; VII}_B$ (0.3.8).

\medskip

For later purposes let us note that any $o$-module $M$ in
$\Mod_\top (o)$ has a unique largest quotient module $M_\cot$
which is Hausdorff and torsionfree: If $(M_j)_{j\in J}$ is the
family of all torsionfree Hausdorff quotient modules of $M$ then
$M_\cot$ is the coimage of the natural map
$M\longrightarrow\prod\limits_{j\in J} M_j$.

For any $o$-module $M$ in $\Mod_{\comp}^{\rm fl} (o)$ we can
construct the $K$-Banach space
$$
M^d:=\Hom^{\cont}_o (M,K)\ \ \hbox{with norm}\ \ \|\ell
\|:=\mathop{\max}\limits_{m\in M} |\ell (m)|\ .
$$
This defines a contravariant additive functor
$$
\matrix{
\Mod^{\rm fl}_{\comp} (o) & \longrightarrow & \Ban (K)\cr\cr
\hfill M & \longmapsto & M^d\hfill}
$$
into the category $\Ban (K)$ of all $K$-Banach spaces with
morphisms being all continuous $K$-linear maps. Actually all maps
in the image of this functor are norm decreasing. The groups of
homomorphisms in $\Mod^{\rm fl}_{\comp} (o)$ are $o$-modules
whereas in $\Ban (K)$ they are $K$-vector spaces. The above
functor therefore extends naturally to a contravariant additive
functor
$$
\Mod^{\rm fl}_\comp (o)_{\Qdss} \;\longrightarrow\; \Ban (K)\ .
$$
Here $\Afr_\Qdss$, for any additive category $\Afr$, denotes the
additive category with the same objects as $\Afr$ and such that
$$
\Hom_{\Afr_\Qdss } (A,B)\; :=\; \Hom_\Afr (A,B)\otimes\Qdss
$$
for any two objects $A,B$ in $\Afr$ with the composition of
morphisms in $\Afr_\Qdss$ being the $\Qdss$-linear extension of
the composition in $\Afr$.

\medskip

{\bf Theorem 1.2:}

{\it The functor
$$
\matrix{
\Mod^{\rm fl}_\comp (o)_\Qdss & \buildrel\sim\over\longrightarrow &
\Ban (K)\cr\cr \hfill M & \longmapsto & M^d \hfill}
$$
is an anti-equivalence of categories.}

Proof: Let $\Ban (K)^{\le 1}$ denote the category of all
$K$-Banach spaces $(E,\|\ \|)$ such that $\|E\|\subseteq |K|$ with
morphisms being all norm decreasing $K$-linear maps. Clearly our
functor factorizes into
$$
\Mod^{\rm fl}_\comp (o) \; \mathop{\longrightarrow}\limits^{(.)^d}\;
\Ban (K)^{\le 1} \; \mathop{\longrightarrow}\limits^{\rm forget}\;
\Ban (K)\ .
$$
For any $K$-Banach space $(E,\|\ \|)$ we may define by $\|v\|' :=
{\rm inf}\{r \in |K|: r \geq \|v\|\}$ another norm $\|\
\|'$ on $E$ satisfying $\|E\|'\subseteq |K|$. Because of
$|\pi | \leq \|v\|/\|v\|' \leq 1$ for $v \neq 0$, where $\pi$ is a
prime element of $K$, the two norms $\|\ \|$ and $\|\ \|'$ are
equivalent. It follows that the right hand functor above induces
an equivalence of categories
$$
(\Ban (K)^{\le 1})_{\Qdss} \buildrel\sim\over\longrightarrow \Ban
(K)\ .
$$
We therefore are reduced to show that
$$
\matrix{
\Mod^{\rm fl}_\comp (o) & \longrightarrow & \Ban (K)^{\le 1}\cr\cr
\hfill M & \longmapsto & M^d\hfill}
$$
is an anti-equivalence of categories. Let $(E,\|\ \|)$ be a
$K$-Banach space and denote by $E^\circ:=\{v \in E:\|v\| \le 1\}$
its unit ball. Then
$$
E^d:=\Hom_o (E^\circ ,o)\;\; \hbox{with the topology of pointwise
convergence}
$$
is a linear-topological $o$-module which is torsionfree and
complete. In fact, $E^d$ is the unit ball of the dual Banach space
$E'$ but equipped with the weak topology. Since
$$
\matrix{
E^d & \hookrightarrow & \prod\limits_{v \in E^\circ }o\cr\cr
\hfill\lambda & \longmapsto & (\lambda (v))_v\hfill}
$$
is a topological embedding we see that $E^d$ is compact. This
defines a functor
$$
\matrix{
\Ban (K)^{\le 1} & \longrightarrow & \Mod^{\rm fl}_\comp (o)\cr\cr
\hfill (E,\|\ \|) & \longmapsto & E^d\ .\hfill}
$$
It is an immediate consequence of Remark 1.1 that, for an
$o$-module $M$ in $\Mod^{\rm fl}_\comp (o)$, the $o$-linear map
$$
\matrix{
\iota_M:M & \longrightarrow & (M^d)_s'\hfill\cr\cr
\hfill m & \longmapsto & [\ell\longmapsto \ell(m)]\hfill}
$$
into the weak dual $(M^d)_s'$ of the Banach space $M^d$ is
injective. Since it is easily seen to be continuous the
compactness of $M$ implies that $\iota_M$ is a closed embedding.
By definition the image of $\iota_M$ is contained in $M^{dd}$.
Assume now that there is a $\lambda\in M^{dd}\setminus {\rm im}
(\iota_M)$. Since ${\rm im} (\iota_M)$ is closed in $(M^d)_s'$
there is, by Hahn-Banach ([Mon] V.1.2 Thm. 5(ii) or [NFA] 13.3), a
continuous linear form on $(M^d)_s'$ which in absolute value is
$\ge 1$ on $\lambda $ and is $<1$ on ${\rm im} (\iota_M)$. But, as
another consequence of Hahn-Banach ([NFA] 9.7), any continuous
linear form on $(M^d)_s'$ is given by evaluation in a vector in
$M^d$. Hence we find an $\ell\in M^d$ such that $|\lambda (\ell
)|\ge 1$ and $|\ell (M)|<1$. The latter implies $\|\ell\|<1$ so
that $\|\lambda (\ell )\|\le \|\lambda
\|\cdot \|\ell\||<1$ which is a contradiction. We obtain that
$\iota_M:M\buildrel\cong\over\longrightarrow M^{dd}$, in fact, is
a topological isomorphism. This means that the $\iota_M$
constitute a natural isomorphism between the identity functor and
the functor $(.)^{dd}$ on $\Mod^{\rm fl}_\comp (o)$. On the other
hand any $(E,\|\ \|)$ in $\Ban (K)^{\le 1}$ is isometric to a
Banach space $c_0 (I)$ for some set $I$ ([Mon] IV.3 Cor. 1 or
[NFA] 10.1). A straightforward explicit computation shows that
$c_0 (I)^{dd}=c_0 (I)$. The functor $M\longrightarrow M^d$
therefore is fully faithful as well as essentially surjective and
consequently an equivalence.

\medskip

The exactness properties of this functor are as follows.

\medskip

{\bf Proposition 1.3:}

{\it For any map $f : M \longrightarrow N$ in $\Mod^{\rm
fl}_{\comp} (o)$ we have:

i. $\ker (f)^d = M^d /\overline{f^d (N^d)}$;

ii. $[\coker(f)_{\cot}]^d = \ker(f^d)$ ;

iii. $f$ is surjective if and only if $f^d$ is an isometry.}

Proof: i. The submodules $\ker(f)$ and $\im(f)$ lie again in
$\Mod^{\rm fl}_{\comp} (o)$. It follows from [SGA3] ${\rm Exp.\;
VII}_B$ (0.3.7) that the surjection $M \longrightarrow \im(f)$
splits, i.e., we have $M \cong \ker(f) \oplus \im(f)$ in
$\Mod^{\rm fl}_{\comp} (o)$. It suffices therefore to consider the
case where $f$ is injective and to show that then the image of
$f^d$ is dense in $M^d$. If not we find by Hahn-Banach a nonzero
continuous linear form $\lambda$ on $M^d$ which vanishes on the
image of $f^d$. Up to scaling we may assume that $\lambda \in
M^{dd}$, i.e., that there is a nonzero $m \in M$ such that
$\lambda(\ell) = \ell(m)$. The vanishing property of $\lambda$
means of course that $f(m) = 0$ which is a contradiction.

iii. If $f$ is surjective then $f^d$ is an isometry by
construction. Suppose now that $f^d$ is an isometry. Let $n \in
N$; we view $n$ as a linear form in the unit ball of the dual
Banach space $(N^d)'$. By Hahn-Banach $n$ extends (via $f^d$) to a
linear form in the unit ball of $(M^d)'$; this means of course
that we find an $m \in M$ such that $f(m) = n$.

ii. Let $E$ denote the kernel of $f^d$. Then $E^d$ is, by iii., a
torsionfree Hausdorff quotient of $\coker(f)$. On the other hand
$[\coker(f)_{\cot}]^d$ clearly is a subspace of $\ker(f^d)$.

\medskip

Let $M$ be a module in $\Mod^{\rm fl}_{\comp} (o)$. Since $M$ is
torsionfree it is an $o$-submodule of the $K$-vector space $M_K := M
\otimes_o K$. Thm. 1.2 tells us that there is a natural identification
$$
M_K = \Hom_o^{\rm cont}(o,M) \otimes \Qdss = \Hom_K^{\rm cont}(M^d,K)
= (M^d)'
$$
between $M_K$ and the continuous dual $(M^d)'$ of the
Banach space $M^d$. We always equip $M_K$ with the finest locally
convex topology such that the inclusion $M
\subseteq M_K$ is continuous. An $o$-submodule $L
\subseteq M_K$ is open if and only if $\alpha L \cap M$ is open in
$M$ for any $0 \neq \alpha \in o$. By construction this topology
has the property that
$$
\Hom^{\cont}_o (M,V) = \L (M_K,V)
$$
for any locally convex $K$-vector space $V$ where, following a
common convention, we let $\L (.,.)$ denote the vector space of
continuous linear maps between two locally convex $K$-vector
spaces. In particular $M^d$ at least as a vector space is the
continuous dual $(M_K)'$. Since under the identification $M_K =
(M^d)'$ the topology of $M$ is induced by the weak topology on
$(M^d)'$ we also see that the identification map $M_K
\mathop{\longrightarrow}\limits^{=} (M^d)_s'$ is continuous. This
shows that $M_K$ is Hausdorff and that $M_K$ also induces the
given topology on $M$.

\medskip

{\bf Lemma 1.4}

{\it The locally convex $K$-vector space $M_K$, for any $M$ in
$\Mod^{\rm fl}_{\comp} (o)$, is complete.}

Proof: Fix a prime element $\pi$ of $K$. Let $\Fscr$ be a minimal
Cauchy filter on $M_K$. We first show that there is a $m\in\Ndss$
such that
$$
F\cap \pi^{-m} M \neq\emptyset\ \ \hbox{for all}\ F\in\Fscr\ .
$$
Otherwise there exists for any $n\in\Ndss$ a $F_n\in\Fscr$ with
$F_n\cap \pi^{-n} M = \emptyset$. By the minimality of $\Fscr$ we
may assume that
$$
F_n=F_n+L_n\ \ \hbox{for some open $o$-submodule}\ L_n\subseteq
M_K
$$
([B-GT] Chap. II \S3.2 Prop. 5). We also may assume that the $L_n$
form a decreasing sequence $L_1\supseteq L_2\supseteq\ldots$. The
$o$-submodule
$$
L := \sum_{n\in\Ndss}(L_n\cap \pi^{-n} M)
$$
is open in $M_K$ since $L\cap \pi^{-n} M \supseteq L_n\cap
\pi^{-n} M$ for all $n\in\Ndss$. The $L_n$ being decreasing and
the $\pi^{-n} M$ being increasing it is clear that
$$
L \subseteq L_n + \pi^{-n} M\ \ \hbox{for all}\ n\in\Ndss\ .
$$
As a Cauchy filter $\Fscr$ must contain a coset $v + L$ for some
$v \in M_K$. If $n_0\in\Ndss$ is chosen in such a way that $v \in
\pi^{-n_0} M$ we have $L_{n_0} + \pi^{-n_0} M \in\Fscr$. Both sets $F_{n_0}$ and
$L_{n_0} + \pi^{-n_0} M$ belonging to the filter $\Fscr$ we obtain
$F_{n_0}\cap (L_{n_0} + \pi^{-n_0} M) \neq \emptyset$, i.e.,
$F_{n_0}\cap \pi^{-n_0} M = (F_{n_0}+L_{n_0}) \cap \pi^{-n_0} M
\neq\emptyset$ which is a contradiction. We see that
$$
\Fscr_m:=\{ F\cap \pi^{-m} M : F\in\Fscr\}
$$
for an appropriate $m\in\Ndss$ is a filter on $\pi^{-m} M$. Since
$\pi^{-m} M$ is compact in $M_K$ the filter $\Fscr_m$ being also a
Cauchy filter has to be convergent. By [B-GT] Chap. II
\S3.2 Cor. 3 then $\Fscr$ is convergent, too. This proves that
$M_K$ is complete.

\medskip

{\bf Lemma 1.5:}

{\it For any two $o$-modules $M$ and $N$ in $\Mod^{\rm fl}_{\comp}
(o)$ we have:

i. For any compact subset $C \subseteq N_K$ the closed
$o$-submodule in $N_K$ topologically generated by $C$ is compact
as well;

ii. for any compact subset $C \subseteq N_K$ there is a $0 \neq
\alpha \in o$ such that $\alpha C \subseteq N$;

iii. $\Hom^{\cont}_o (M,N) \otimes \Qdss = \L (M_K,N_K)$;

iv. passing to the transpose induces a $K$-linear isomorphism}
$$
\L (M_K,N_K) \mathop{\longrightarrow}\limits^{\cong} \L (N^d,M^d)\
.
$$

Proof: i. Let $<C>$ denote the $o$-submodule generated by $C$. Let
$L\subseteq N_K$ be any open and therefore also closed
$o$-submodule. Since $C$ is compact we find finitely many
$c_1,\ldots,c_n\in C$ with $C\subseteq (c_1+L)\cup\ldots\cup
(c_n+L)$. Then $<C>$ is contained in $oc_1 + \ldots + oc_n + L$.
But $oc_1 + \ldots + oc_n$ is compact, too, so that we again find
finitely many $a_1,\ldots,a_m\in oc_1 + \ldots + oc_n$ with
$$
oc_1 + \ldots + oc_n \subseteq (a_1+L)\cup\ldots\cup (a_m+L)\ .
$$
Together we obtain
$$
<C>\ \subseteq \mathop{\bigcup}\limits_{1\leq i\leq m}a_i+L
$$
and since the right hand side is closed the closure
$\overline{<C>}$ of $<C>$ also satisfies
$$
\overline{<C>}\ \subseteq \bigcup\limits_{1\leq i\leq m}a_i+L\ .
$$
Since $L$ was arbitrary this implies by [B-GT] Chap. II \S4.2 Thm.
3 that $\overline{<C>}$ is precompact. On the other hand, as a
consequence of Lemma 1.4, $\overline{<C>}$ is Hausdorff and
complete. Hence $\overline{<C>}$ is compact.

ii. By i. we may assume that $C$ is a compact $o$-submodule of
$N_K$. Fix a prime element $\pi$ of $K$ and put $C_n := C \cap
\pi^{-n}N$ for any $n \in \Ndss$. These $C_n$ form an increasing
sequence $C_1\subseteq C_2\subseteq\ldots$ of compact
$o$-submodules of $C$ such that $C=\bigcup\limits_{n\in\Ndss}C_n$.
We have to show that $C_m=C$ holds for some $m\in\Ndss$. Being
empty the subset $\bigcap\limits_{n\in\Ndss}(C\backslash C_n)$ is
not dense in $C$. As a compact space $C$ in particular is a Baire
space ([B-GT] Chap. IX \S5.3 Thm. 1) so that already some
$C\backslash C_n$ is not dense in $C$. This means that $C_n$
contains a non-empty open subset of $C$. It is then itself an open
$o$-submodule and therefore has to be of finite index in $C$. Our
claim obviously follows from that.

iii. We have $\Hom^{\cont}_o (M,N) \otimes \Qdss = \Hom^{\cont}_o
(M,N) \otimes_o K = \Hom^{\cont}_o (M,N_K) = \L (M_K,N_K)$ where
the second identity is a consequence of the second assertion.

iv. This follows from iii. and Thm. 1.2.

\medskip

The assertion iii. in Lemma 1.5 in particular means that $M_K$ and
$N_K$ are isomorphic in $\Mod^{\rm fl}_\comp (o)_\Qdss$ if and
only if they are isomorphic as locally convex vector spaces.

For any two $o$-modules $M$ and $N$ in $\Mod^{\rm fl}_\comp (o)$
we always view $\Hom^{\cont}_o (M,N)$ as a linear-topological
$o$-module by equipping it with the topology of compact
convergence. As a consequence of Lemma 1.5 this topology is
induced by the topology of compact convergence on the vector space
$\L(M_K,N_K)$. We write $\L_{cc}(M_K,N_K)$ for $\L(M_K,N_K)$
equipped with the finest locally convex topology such that the
inclusion $\Hom^{\cont}_o (M,N) \hookrightarrow \L_{cc}(M_K,N_K)$
is continuous. By a similar argument as before $\L_{cc}(M_K,N_K)$
is Hausdorff and the latter inclusion is a topological embedding.
Moreover, by Lemma 1.5, $\L_{cc}(M_K,N_K)$ is, in both variables,
a functor on $\Mod^{\rm fl}_\comp (o)_\Qdss$.

Given two $K$-Banach spaces $E_1$ and $E_2$ we write, following
traditional usage, $\L_{s}(E_1,E_2)$ for the vector space
$\L(E_1,E_2)$ equipped with the locally convex topology of
pointwise convergence. We write $\L_{bs}(E_1,E_2)$ for
$\L(E_1,E_2)$ equipped with the finest locally convex topology
which coincides with the topology of pointwise convergence on any
equicontinuous subset in $\L(E_1,E_2)$. Corresponding to any
choice of defining norms $\|\ \|_i$ on $E_i$ for $i = 1,2$ we have
the operator norm $\|\ \|$ on $\L(E_1,E_2)$. A subset in
$\L(E_1,E_2)$ is equicontinuous if and only if it is bounded with
respect to $\|\ \|$. Hence the topology of $\L_{bs}(E_1,E_2)$ can
equivalently be characterized as being the finest locally convex
topology which induces the topology of pointwise convergence on
the unit ball with respect to $\|\ \|$ in $\L(E_1,E_2)$.

\medskip

{\bf Proposition 1.6:}

{\it Passing to the transpose induces, for any $o$-modules $M$ and
$N$ in $\Mod^{\rm fl}_\comp (o)$, an isomorphism of locally convex
$K$-vector spaces}
$$
\L_{cc}(M_K,N_K) \mathop{\longrightarrow}\limits^{\cong}
\L_{bs}(N^d,M^d)\ .
$$

Proof: It is clear from our preliminary discussion that it
suffices to show that the $o$-linear isomorphism
$$
\Hom^{\cont}_o (M,N)\;\mathop{\longrightarrow}\limits^{\cong}\;
\{f \in \L(N^d,M^d) : \|f\|\leq 1\}
$$
given by the transpose is topological provided the left, resp.
right, hand side carries the topology of compact, resp. pointwise
convergence. We recall from the proof of Thm. 1.2 that $M$ is the
unit ball in the dual Banach space $(M^d)'$ equipped with the weak
topology; we also have seen there that the closed equicontinuous
subsets of the weak dual $(M^d)_s'$ are compact. By the
Banach-Steinhaus theorem ([Tie] Thm. 4.3) a subset of $(M^d)_s'$
is equicontinuous if and only if it is bounded. Clearly any
compact subset is bounded. It follows that for a closed subset of
$(M^d)_s'$ the following properties are equivalent: Bounded,
equicontinuous, bounded for the dual Banach norm, compact. This
shows that the topology of compact convergence on $\Hom^{\cont}_o
(M,N)$ is induced by the strong topology on
$\L((M^d)_s',(N^d)_s')$. Our assertion therefore will be a
consequence of the quite general fact that for any two $K$-Banach
spaces $E_1$ and $E_2$ the transpose induces a topological
isomorphism
$$
\L_s(E_1,E_2) \mathop{\longrightarrow}\limits^{\cong}
\L_b((E_2)_s',(E_1)_s')
$$
where on the right hand side the subscript $b$ indicates, as
usual, the strong topology. This is straightforward from the
definitions and the fact that
$$
\matrix{\hbox{set of all open}\hfill &
\mathop{\longrightarrow}\limits^{\sim} & \hbox{set of all
closed equicontinuous}\hfill\cr o\hbox{-submodules in}\ E_2 & &
o\hbox{-submodules in}\ (E_2)_s'\hfill \cr\cr
\hfill L & \longmapsto & L^p := \{\ell \in (E_2)' : |\ell(v)|\leq 1\
\hbox{for any}\ v \in L\}\hfill }
$$
is a bijection which is a direct consequence of the Hahn-Banach
theorem.

\bigskip

{\bf 2. Iwasawa modules and representations}

\smallskip

>From now on we let $G$ denote a fixed profinite group. The
completed group ring of $G$ (over $o$) is defined to be
$$
o[[G]] := \mathop{\lim\limits_{\longleftarrow}}\limits_{H
\in \Nscr} o[G/H]
$$
where $\Nscr = \Nscr (G)$ denotes the family of all open normal
subgroups of $G$. In a natural way $o[[G]]$ is a torsionfree and
compact linear-topological $o$-module; the ring multiplication is
continuous. The surjections $o[G] \longrightarrow o[G/H]$ for $H
\in \Nscr$ induce in the limit a ring homomorphism
$$
o[G] \longrightarrow o[[G]]
$$
whose image is dense and which is injective ([Laz] II.2.2.3.1).
Being the projective limit of the inclusions $G/H \subseteq
o[G/H]$ the composed map
$$
G \mathop{\longrightarrow}\limits^{\subseteq} o[G] \longrightarrow
o[[G]]
$$
is continuous and hence, by compactness, a homeomorphism onto its
image.

Consider now a module $M$ in $\Mod_{\top} (o)$ and let $C(G,M)$
denote the $o$-module of all continuous maps from $G$ into $M$. It
follows from the above discussion that the $o$-linear map
$$
\matrix{
\Hom^{\cont}_o (o[[G]],M) & \longrightarrow & C(G,M)\cr\cr
\hfill f & \longmapsto & f|G \hfill}
$$
is well defined and injective.

\eject

{\bf Lemma 2.1:}

{\it For any complete $o$-module $M$ in $\Mod_{\top} (o)$ the map
$$
\matrix{
\Hom^{\cont}_o (o[[G]],M) & \mathop{\longrightarrow}\limits^{\cong}
 & C(G,M)\cr\cr \hfill f & \longmapsto & f|G \hfill}
$$
is a bijection.}

Proof: We extend a given $\varphi \in C(G,M)$ $o$-linearly to
$o[G]$. By the completeness assumption and the density of $o[G]$
in $o[[G]]$ it suffices to show that this extension, which we
again denote by $\varphi$, is continuous with respect to the
topology induced by $o[[G]]$. Fix an open $o$-submodule $L
\subseteq M$. By the uniform continuity of $\varphi$ on $G$ we
find an $H \in
\Nscr$ such that
$$
\varphi(g_iH) \subseteq \varphi(g_i)+L
$$
for all $g_i$ in a system of representatives for the left cosets
of $H$ in $G$ (compare [Laz] II.2.2.5). Let $\alpha \in o$ be some
element such that $\alpha\cdot \varphi(g_i)\subseteq L$ for all
$i$. The $o$-submodule
$$
L':=\mathop{\oplus}\limits_{i}\,\{ \sum_{g\in g_iH} r_gg :\sum_g
r_g\in \alpha o\}
$$
then is open in $o[G]$ and we have
$$
\varphi(L')\subseteq\sum_i(\alpha o\cdot \varphi(g_i)+L)\subseteq L\ .
$$

\medskip

We set $K[[G]] := o[[G]]_K$. This is a locally convex vector space
as well as a $K$-algebra such that the multiplication is
separately continuous.

\medskip

{\bf Corollary 2.2}

{\it For any quasi-complete Hausdorff locally convex $K$-vector
space $V$ we have the $K$-linear isomorphism}
$$
\matrix{
\L (K[[G]],V) & \mathop{\longrightarrow}\limits^{\cong}
 & C(G,V)\cr\cr \hfill f & \longmapsto & f|G \hfill\ .}
$$

Proof: The map is clearly well defined and injective. For the
surjectivity let $\varphi \in C(G,V)$. Define $M$ to be the closed
$o$-submodule of $V$ topologically generated by $\varphi (G)$.
This $M$ lies in $\Mod_{\top} (o)$. Since $G$ is compact $M$ is
bounded in $V$. The quasi-completeness of $V$ therefore ensures
that $M$ is complete. Hence we have, by Lemma 2.1, a continuous
$o$-linear map $f : o[[G]] \longrightarrow M \subseteq V$ such
that $f|G = \varphi$. The $K$-linear extension of $f$ then is the
preimage of $\varphi$ we were looking for.

\medskip

We will apply these results to obtain a $G$-equivariant version of
the duality theorem of the previous section.

\medskip

{\bf Definition:}

{\it A $K$-Banach space representation $E$ of $G$ is a $K$-Banach
space $E$ together with a $G$-action by continuous linear
automorphisms such that the map $G \times E \longrightarrow E$
describing the action is continuous.}

\medskip

We define
$$
\matrix{
\Ban_G(K) := & \hbox{category of all}\ K\hbox{-Banach
representations of}\ G\ \hbox{with}\hfill\cr & \hbox{morphisms
being all}\ G\hbox{-equivariant continuous linear maps.}\hfill}
$$
As a consequence of the Banach-Steinhaus theorem ([Tie] Thm.
4.1.1$^{\circ}$), to give a $K$-Banach representation of $G$ on
the $K$-Banach space $E$ is the same as to give a continuous
homomorphism $G \longrightarrow \L_s(E,E)$. But $\L_s(E,E)$ is
quasi-complete and Hausdorff ([B-TVS] III.27 Cor. 4 or [NFA]
7.14). Hence such a homomorphism extends, by Cor. 2.2, uniquely to
a continuous $K$-linear map $K[[G]]
\longrightarrow \L_s(E,E)$. By a density argument the
latter map is a $K$-algebra homomorphism. This shows that a
$K$-Banach space representation of $G$ on $E$ is the same
as a separately continuous action $K[[G]] \times E
\longrightarrow E$ of the algebra $K[[G]]$ on $E$.

Since the image of $o[[G]]$ in $\L_s(E,E)$ under the above
homomorphism is compact and hence (by Banach-Steinhaus)
equicontinuous we also have that a $K$-Banach space representation
of $G$ on $E$ is the same as a continuous (unital) homomorphism of
$K$-algebras $K[[G]] \longrightarrow \L_{bs}(E,E)$.

\medskip

{\bf Definition:}

{\it An Iwasawa $G$-module over $o$ is an $o$-module $M$ in
$\Mod^{\rm fl}_\comp (o)$ together with a continuous (left) action
$o[[G]] \times M \longrightarrow M$ of the compact $o$-algebra
$o[[G]]$ on $M$ such that the induced $o$-action on $M$ is the
given $o$-module structure.}

\medskip

Let
$$
\matrix
{\Mod^{\rm fl}_\comp (o[[G]]) := & \hbox{category of all Iwasawa}\
G\hbox{-modules over}\ o\hfill\cr & \hbox{with morphisms being all
continuous}\ o[[G]]\hbox{-}\hfill\cr & \hbox{module
homomorphisms.}\hfill }
$$
A continuous (unital) homomorphism of $K$-algebras
$$
K[[G]] \longrightarrow \L_{cc}(M_K,M_K)\leqno{(\ast)}
$$
for some $M_K$ in $\Mod^{\rm fl}_\comp (o)_{\Qdss}$ induces a
continuous map $o[[G]] \longrightarrow \Hom^{\cont}_o (M,\break
M_K)$ where the right hand side carries the topology of compact
convergence. By [B-GT] Chap. X \S3.4 Thm. 3 this is the same as a
continuous map $o[[G]] \times M \longrightarrow M_K$. According to
Lemma 1.5.ii the image of this latter map is contained in
$\alpha^{-1}M$ for some $0 \neq \alpha \in o$. If $N$ denotes the
closed $o$-submodule of $M_K$ topologically generated by this
image we therefore have $N_K = M_K$, and the above homomorphism of
$K$-algebras $(\ast)$ is the tensor product with $K$ of a
continuous (unital) homomorphism of $o$-algebras
$$
o[[G]] \longrightarrow \Hom^{\cont}_o (N,N)\ .
$$
Again by [B-GT] loc. cit. this is the same as a continuous action
$o[[G]] \times N \longrightarrow N$ of the compact $o$-algebra
$o[[G]]$ on the $o$-module $N$ in $\Mod^{\rm fl}_\comp (o)$ which
extends the $o$-module structure.

Hence we see that to give a continuous (unital) homomorphism of
$K$-algebras $(\ast)$ is the same as to give an object in the
category $\Mod^{\rm fl}_\comp (o[[G]])_{\Qdss}$. By combining this
discussion with Prop. 1.6 we arrive at the following equivariant
version of Thm. 1.2.

\medskip

{\bf Theorem 2.3:}

{\it The functor
$$
\matrix{
\Mod^{\rm fl}_\comp (o[[G]])_\Qdss & \buildrel\sim\over\longrightarrow &
\Ban_G (K)\cr\cr \hfill M & \longmapsto & M^d \hfill}
$$
is an anti-equivalence of categories.}

\bigskip

{\bf 3. Admissible representations}

\smallskip

In order to obtain a reasonable theory of Banach space
representations it seems necessary to impose certain additional
finiteness conditions. The first idea is to consider only those
$K$-Banach space representations of $G$ which correspond, under
the duality of the previous section, to finitely generated
$K[[G]]$-modules. As a consequence of the compactness of the ring
$o[[G]]$ it will turn out that the theory of these representations
in fact is completely algebraic in nature. In order to obtain an
intrinsic characterization we will assume in this section that $G$
is a compact $p$-adic Lie group. We then have:

- The subfamily of all topologically finitely generated pro-$p$-groups
in $\Nscr = \Nscr(G)$ is cofinal ([B-GAL] Chap. III \S1.1 Prop. 2(iii)
and \S7.3 and 4 and [Laz] III.2.2.6 and III.3.1.3).

- The ring $o[[G]]$ is left and right noetherian ([Laz] V.2.2.4).

The ring $K[[G]]$ then is left and right noetherian as well.

\medskip

{\bf Definition:}

{\it A $K$-Banach space representation $E$ of $G$ is called admissible
if there is a $G$-invariant bounded open $o$-submodule $L \subseteq E$
such that, for any $H \in \Nscr$, the $o$-submodule $(E/L)^H$ of
$H$-invariant elements in the quotient $E/L$ is of cofinite type.}

\medskip

We recall that an $o$-module $N$ is called of cofinite type if its
Pontrjagin dual $\Hom_{\rm o}(N,K/o)$ is a finitely generated
$o$-module. We also point out that an arbitrary open $o$-submodule
$L \subseteq E$ contains the $G$-invariant open $o$-submodule
$\mathop{\bigcap}\limits_{g \in G} gL$.

Let
$$
\matrix{
\Ban_G^{\rm adm}(K) := & \hbox{the full subcategory in}\
\Ban_G(K)\hfill\cr & \hbox{of all admissible representations.
\hfill} }
$$
On the other hand we let $\Mod_{\rm fg}^{\rm fl}(o[[G]])$, resp.
$\Mod_{\rm fg}(K[[G]])$, denote the category of all finitely
generated and $o$-torsionfree (left unital) $o[[G]])$-modules,
resp. of all finitely generated (left unital) $K[[G]])$-modules.
It is clear that
$$
\Mod_{\rm fg}^{\rm fl}(o[[G]])_{\Qdss} =
\Mod_{\rm fg}(K[[G]])\ .
$$
Since $K[[G]]$ is noetherian the category $\Mod_{\rm fg}(K[[G]])$ is
abelian.

\medskip

{\bf Proposition 3.1:}

{\it i. A finitely generated $o[[G]]$-module $M$ carries a unique
Hausdorff topology - its canonical topology - such that the action
$o[[G]] \times M \longrightarrow M$ is continuous;

ii. any submodule of a finitely generated $o[[G]]$-module is
closed in the canonical topology;

iii. any $o[[G]]$-linear map between two finitely generated
$o[[G]]$-modules is continuous for the canonical topologies.}

Proof: Since $o[[G]]$ is compact and noetherian this is an easy
exercise. But we point out that the assertions hold for any
compact ring by [AU] Cor. 1.10.

\medskip

It follows that equipping a module in $\Mod_{\rm fg}^{\rm
fl}(o[[G]])$ with its canonical topology induces a fully faithful
embedding $\Mod_{\rm fg}^{\rm fl}(o[[G]]) \longrightarrow
\Mod^{\rm fl}_\comp (o[[G]])$. This then in turn induces a fully
faithful embedding $\Mod_{\rm fg}(K[[G]]) \longrightarrow
\Mod^{\rm fl}_\comp (o[[G]])_{\Qdss}$. In other words we can and
will view $\Mod_{\rm fg}(K[[G]])$ as a full subcategory of
$\Mod^{\rm fl}_\comp (o[[G]])_{\Qdss}$.

For each $H \in \Nscr$ let $I_H$ denote the kernel of the
projection map $o[[G]] \longrightarrow o[G/H]$. This is a family
of 2-sided ideals in $o[[G]]$ which converges to zero. As a left
(or right) ideal $I_H$ is generated by the elements $h - 1$ for $h
\in H$. For the sake of completeness we include a proof of the following
well known fact.

\medskip

{\bf Lemma 3.2:}

{\it Let $H \in \Nscr$ be a pro-$p$-group; then the ideal powers
$I_{H}^n$, for $n \in \Ndss$, converge to zero.}

Proof: We may assume that $G$ is finite. Let $\pi$ denote a prime
element in $o$ and $k := o/\pi o$ the residue field of $o$. By
Clifford's theorem ([CR] (49.2)) and [Ser] IX\S1 the ideal ${\rm
ker}(k[G] \longrightarrow k[G/H])$ is contained in the radical of
the ring $k[G]$. Since this radical is nilpotent we have $I_{H}^m
\subseteq \pi o[G]$ for some $m \in \Ndss$.

\medskip

{\bf Lemma 3.3:}

{\it Let $H \in \Nscr$ be a pro-$p$-group; a module $M$ in
$\Mod^{\rm fl}_\comp (o[[G]])$ is finitely generated over $o[[G]]$
if and only if $M/I_{H} M$ is finitely generated over $o$.}

Proof: This is the well known Nakayama lemma; compare [BH] for a
thorough discussion.

\medskip

{\bf Lemma 3.4:}

{\it A $K$-Banach space representation $E$ of $G$ is admissible if
and only if the dual space $E'$ is finitely generated over
$K[[G]]$.}

Proof: Let us first assume that $E'$ is finitely generated over
$K[[G]]$. There is then a finitely generated $o[[G]]$-submodule $M
\subseteq E'$ such that $E' = M_K$. After equipping $M$ with its
canonical topology we have $E = M^d$. Moreover $L :=
\Hom^{\cont}_o (M,o)$ is a $G$-invariant bounded open $o$-submodule
in $E$. It follows from Remark 1.1 that $E/L = \Hom^{\cont}_o
(M,K/o)$ (where $K/o$ carries the discrete topology) and hence
that
$$
(E/L)^H = \Hom^{\cont}_o (M,K/o)^H = \Hom^{\cont}_o (M/I_H
M,K/o)\leqno{(\ast)}
$$
for any $H \in \Nscr$. Hence $(E/L)^H$ is of cofinite type.

On the other hand let now $H \in \Nscr$ be a pro-$p$-group and $L
\subseteq E$ be a $G$-invariant bounded open $o$-submodule such that
$(E/L)^H$ is of cofinite type. In the proof of Prop. 1.6 we had
recalled that the $G$-invariant $o$-submodule $M := L^p$ in $E_s'$ is
compact. Since $L$ is bounded we have $E' = M_K$. So the identities
$(\ast)$ apply correspondingly and we obtain that $\Hom^{\cont}_o (
M/I_H M,K/o)$ is of cofinite type. But since $I_H$ is finitely
generated as a right ideal the submodule $I_H M$ is the image of
finitely many copies $M \times \ldots \times M$ under a continuous map
and hence is closed in $M$. By Pontrjagin duality and the Nakayama
lemma over $o$ applied to the compact $o$-module $M/I_H M$ the latter
therefore is finitely generated over $o$. Lemma 3.3 then implies that
$M$ is finitely generated over $o[[G]]$ and hence that $E'$ is
finitely generated over $K[[G]]$.

\medskip

The above proof shows that the defining condition for admissibility
only needs to be tested for a single pro-$p$-group $H \in \Nscr$. On
the other hand assume $E$ to be an admissible representation of $G$
and let $L \subseteq E$ be as in the above definition. Consider an $H
\in \Nscr$ and an {\it arbitrary} $G$-invariant open $o$-submodule
$L_{\rm o} \subseteq E$. We claim that $(E/L_{\rm o})^H$ is of
cofinite type. Replacing $L$ by $\alpha L$ for some appropriate $0
\neq \alpha \in o$ we may assume that $L \subseteq L_{\rm o}$. As we
have seen in the above proof $M := L^p$ is a finitely generated
$o[[G]]$-module. Since $o[[G]]$ is noetherian the $o[[G]]$-submodule
$M_{\rm o} := L_{\rm o}^p$ of $M$ also is finitely generated. As we
have seen this implies that $(E/L_{\rm o})^H$ is of cofinite type.

\medskip

{\bf Theorem 3.5:}

{\it The functor
$$
\matrix{
\Mod_{\rm fg}(K[[G]]) & \buildrel\sim\over\longrightarrow &
\Ban_G^{\rm adm} (K)\cr\cr \hfill M & \longmapsto & M^d \hfill}
$$
is an anti-equivalence of categories.}

Proof: Since $\Mod_{\rm fg}(K[[G]])$ is a full subcategory of
$\Mod^{\rm fl}_\comp (o[[G]])_{\Qdss}$ by Prop. 3.1 this follows from
Thm. 2.3 and Lemma 3.4.

\medskip

As an immediate consequence we obtain that the category
$\Ban_G^{\rm adm} (K)$ is abelian.

\medskip

{\bf Corollary 3.6:}

{\it The functor $E \longmapsto E'$ induces a bijection}
$$
\matrix{
\hbox{set of isomorphism classes} & & \cr \hbox{of topologically
irreducible} & \buildrel\sim\over\longrightarrow & \hbox{set of
isomorphism classes}\cr \hbox{admissible}\ K\hbox{-Banach space} & &
\hbox{of simple}\ K[[G]]\hbox{-modules.}\cr \hbox{representations of\ }G
}
$$

Proof: For any proper closed $G$-invariant subspace $\{0\} \neq
E_{\rm o} \subsetneqq E$ we have, by Hahn-Banach, the exact
sequence of dual vector spaces $0 \rightarrow (E/E_{\rm o})'
\rightarrow E' \rightarrow E_{\rm o}' \rightarrow 0$ in which all
three terms are nonzero. If the $K[[G]]$-module $E'$ is simple the
representation $E$ therefore must be topologically irreducible. On
the other hand write $E' = M_K$ for some module $M$ in $\Mod^{\rm
fl}_{\rm fg} (o[[G]])$ and let $\{0\} \neq V \subsetneqq M_K$ be a
proper $K[[G]]$-submodule. By Prop. 3.1.ii the nonzero
$o[[G]]$-submodule $N := V \cap M$ is closed in $M$ and hence lies
in $\Mod^{\rm fl}_\comp (o[[G]])$. Since the quotient $(M/N)_{\rm
cot} = M/N$ is nonzero as well it follows from Prop. 1.3 that the
kernel of the dual map $E = M^d \longrightarrow N^d$ is a nonzero
proper closed $G$-invariant subspace of $E$.

\medskip

One of the typical pathologies of general Banach space representations
of $G$ is avoided by the admissibility requirement as the following
result shows.

\medskip

{\bf Corollary 3.7:}

{\it Any nonzero $G$-equivariant continuous linear map between two
topologically irreducible admissible $K$-Banach space representations
of $G$ is an isomorphism.}

Proof: This is immediate from Thm. 3.5 and Cor. 3.6.

\medskip

The simplest group to which the results of this section apply is the
group $G = \Zdss_p$ of $p$-adic integers. As is shown in [Dia] already
this group has an extreme wealth of topologically irreducible
$K$-Banach space representations. On the other hand for a commutative
group all ''reasonable'' topologically irreducible $K$-Banach space
representations should be finite dimensional. This is achieved by the
admissibility requirement. The ring $o[[\Zdss_p]]$ is the ring
considered in classical Iwasawa theory; it is isomorphic to the power
series ring $o[[T]]$ in one variable over $o$ ([Was] 7.1). It follows
([Was] \S13.2) that $K[[G]]$ is a principal ideal domain in which
every maximal ideal is of finite codimension.

\medskip

{\bf Remark:} In [ST2] we have introduced the notion of an
analytic module over the algebra $D(G,K)$ of $K$-valued
distributions on $G$ and we have advocated the conjecture that any
$D(G,K)$-module of finite presentation is analytic. Since $K[[G]]$
is naturally a subalgebra of $D(G,K)$ base change would (assuming
this conjecture) induce a functor from $\Mod_{\rm fg}(K[[G]])$
into the category of analytic $D(G,K)$-modules. Since the latter
are dual to a certain class of locally analytic
$G$-representations this functor should correspond to the passage
from a $K$-Banach space representation to the subspace of locally
analytic vectors. The next basic question in this context then
would be whether the ring extension $K[[G]] \longrightarrow
D(G,K)$ is faithfully flat. This is in the spirit of whether every
admissible $K$-Banach space representation of $G$ contains a
locally analytic vector.

\eject

{\bf 4. The group} $G = GL_2 (\dZ)$

\smallskip

In this section we will analyze a certain infinite series of Iwasawa
modules for the group $G := GL_2 (\dZ)$. Let $B \subseteq G$ denote
the Iwahori subgroup of all matrices which are lower triangular modulo
$p$. In $B$ we consider the subgroups $P, P^-$, and $T$ of lower
triangular, upper triangular, and diagonal matrices, respectively. We
also need the subgroups $U$ and $U^-$ of unipotent matrices in $P$ and
$P^-$, respectively. We fix a continuous character $\chi : T
\longrightarrow o^{\times}$. By Cor. 2.2 it extends uniquely to
a continuous homomorphism of $K$-algebras $\chi : K[[T]]
\longrightarrow K$. The inclusions $P \subseteq B \subseteq G$,
resp. the projection $P \longrightarrow T$, induce continuous
algebra monomorphisms $K[[P]] \subseteq K[[B]] \subseteq K[[G]]$,
resp. a continuous algebra epimorphism $K[[P]]
\longrightarrow K[[T]]$. We denote by $\Kchi$ the one
dimensional $K[[P]]$-module given by the composed homomorphism
$K[[P]] \longrightarrow K[[T]]
\mathop{\longrightarrow}\limits^{\chi} K$. Our aim is to study
the finitely generated $K[[G]]$- and $K[[B]]$-modules
$$
M_{\chi} := K[[G]] \mathop{\otimes}\limits_{K[[P]]} \Kchi
\ \ \ \hbox{and}\ \ \ N_{\chi} := K[[B]] \mathop{\otimes}\limits_{K[[P]]}
\Kchi\ ,
$$
respectively. In a similar way (and by a slight abuse of notation)
we have the finitely generated $K[[B]]$-module
$$
N_{\chi}^- := K[[B]] \mathop{\otimes}\limits_{K[[P^-]]}
\Kchi\ .
$$
Put $w := \pmatrix{0 & -1\cr 1 & 0} \in G$ and $w\chi (t) := \chi
(w^{-1}tw)$. As a consequence of the Bruhat decomposition $G = B\
\dot{\cup}\ BwP$ the module $M_{\chi}$, as a $K[[B]]$-module,
decomposes into
$$
M_{\chi} \cong N_{\chi} \oplus N_{w\chi}^-\ .
$$
For later use we note that this decomposition is {\it not}
$K[[G]]$-equivariant since obviously $wN_{\chi} \subseteq
N_{w\chi}^-$.

The module theoretic properties of the series of modules
$N_{\chi}$ and $M_{\chi}$ are governed by one numerical invariant
$c(\chi ) \in K$ of the character $\chi$ which is defined by the
expansion
$$
\chi(\pmatrix{a^{-1} & 0\cr 0 & a})=\exp(c(\chi)\log(a))
$$
for $a$ sufficiently close to 1 (the existence follows from the
topological cyclicity of the group $1+p\dZ$).

In order to investigate the module $N_{\chi}$ we use the Iwahori
decomposition which says that multiplication induces a
homeomorphism $U^- \times\,P
\mathop{\longrightarrow}\limits^{\sim} B$. It implies that $o[[B]]
= o[[U^- ]]\ \widehat{\otimes}\ o[[P]]$ where $\widehat{\otimes}$
is the completed tensor product for linear-topological $o$-modules
([SGA3] ${\rm Exp.\; VII}_B$ (0.3)). The inclusion $K[[U^- ]]
\subseteq K[[B]]$ therefore induces an isomorphism of $K[[U^-
]]$-modules
$$
K[[U^- ]]\;\mathop{\longrightarrow}\limits^{\cong}\;N_{\chi}\
.\leqno{(\ast)}
$$
In particular, any $K[[B]]$-submodule of $N_{\chi}$ corresponds to a
certain ideal in the ring $K[[U^- ]]$. Since the matrix $\gamma :=
\pmatrix{1 & p\cr 0 & 1}$ is a topological generator of $U^-$ the ring
$K[[U^- ]]$ is the ring of formal power series in $\gamma - 1$ whose
coefficients are bounded. As already recalled earlier this is a
principal ideal domain and each ideal is generated by a polynomial in
$\gamma - 1$ with all its zeros lying in the open unit disk ([Was]
\S7.1).

\medskip

{\bf Proposition 4.1:}

{\it If $c(\chi) \not\in \Ndss_{0}$ then $N_{\chi}$ is a simple
$K[[B]]$-module.}

Proof: Let $N \subseteq N_{\chi}$ be a nonzero $K[[B]]$-submodule, $I
\subseteq K[[U^- ]]$ be the ideal which corresponds to $N$ under the
above isomorphism $(\ast)$, and $F_I (\gamma - 1)$ be a polynomial
which generates $I$ and has all its zeros in the open unit disk. The
action of the element $t_a :=
\pmatrix{a & 0\cr 0 & 1} \in T$ on $N_{\chi}$ is given on the left
hand side of $(\ast)$ by
$$
F(\gamma - 1)\ \longmapsto\ \chi(t_a)\cdot F(\gamma^a - 1)
$$
for any bounded power series $F(x)$. Using the bounded power
series
$$
\omega_a(x) := (x + 1)^a - 1 = \sum\limits_{n\in \Ndss} {a\choose n} x^n
$$
this can be rewritten as
$$
F(\gamma - 1)\ \longmapsto\ \chi(t_a)\cdot F(\omega_a(\gamma -
1))\ .
$$
Since this action preserves the ideal $I$ it follows that with $z$
every $\omega_a(z)$, for $a \in \dZ^{\times}$, is a zero of the
polynomial $F_I(x)$. This is only possible if $z+1$ is a $p^m$-th
root of unity for some $m \in \Ndss$. We therefore see that there
are natural numbers $k_{\rm o}$ and $\ell$ such that $F_I(x)$
divides $\omega_{p^{k_{\rm o}}}(x)^{\ell}$. In particular, for any
natural number $k \geq k_{\rm o}$, the polynomial
$\omega_{p^k}(x)^{\ell}$ lies in $I$. We now look at the action of
the element $u := \pmatrix{1 & 0\cr 1 & 1}$ on $N_{\chi}$. It is
straightforward to check that on the left hand side of $(\ast)$ we
have
$$
u(\gamma^n) = \chi(\pmatrix{(1+np)^{-1}& 0\cr 0 &
1+np})\cdot\gamma^{n/(1+np)}\ \ \ \hbox{for any}\ n \in\Ndss_{0}.
$$
It follows that, for $k \geq k_{\rm o}$, with
$\omega_{p^k}(x)^{\ell}$ also
$$
u((\gamma^{p^k} - 1)^{\ell}) = \sum_{j=0}^{\ell} (-1)^j
\pmatrix{\ell\cr j}
\chi(\pmatrix{(1+jp^{k+1})^{-1}& 0\cr 0 & 1+jp^{k+1}})\cdot
\gamma^{jp^k/(1+jp^{k+1})}
$$
lies in the ideal $I$. If $\omega_{p^k}(x)^{\ell}$ and its image
under $u$, for some $k \geq k_{\rm o}$, have no zero in common
then $I$ has to be the unit ideal which means that $N = N_{\chi}$.
In the opposite case we obtain
$$
\matrix{
0 & = & \mathop{\sum}\limits_{j=0}^{\ell} (-1)^j \pmatrix{\ell\cr
j}\chi(\pmatrix{(1+jp^{k+1})^{-1}& 0\cr 0 & 1+jp^{k+1}})\cr\cr\cr
& = & \mathop{\sum}\limits_{j=0}^{\ell} (-1)^j \pmatrix{\ell\cr
j}\exp(c(\chi)\log(1+jp^{k+1}))\hfill }
$$
for any sufficiently big $k \geq k_{\rm o}$. This implies that the
function
$$
\sum_{j=0}^{\ell} (-1)^j \pmatrix{\ell\cr j}\exp(c(\chi)
\log(1+jy))
$$
of the variable $y$ which is analytic in a sufficiently small open
disk around zero has infinitely many zeros and hence vanishes
identically. In order to prove our assertion we therefore have to
show that this is only possible if $c(\chi) \in \Ndss_{0}$. But if
$c(\chi) \not\in \Ndss_{0}$ then evaluating all higher derivatives
of the above function in zero would lead to the identities
$$
\sum_{j=1}^{\ell} (-1)^j \pmatrix{\ell\cr j}j^m = 0\ \ \ \hbox{for
any}\ m \in \Ndss.
$$
This is clearly impossible.

\medskip

The proof of the following companion result is completely
analogous and is therefore omitted.

\medskip

{\bf Proposition 4.2:}

{\it If $c(\chi) \not\in -\Ndss_{0}$ then $N_{\chi}^-$ is a simple
$K[[B]]$-module.}

\medskip

{\bf Lemma 4.3:}

{\it $\Hom_{K[[B]]}(N_{\chi'},N_{\chi}^-) =
\Hom_{K[[B]]}(N_{\chi}^-,N_{\chi'}) = 0$ for any two continuous
characters $\chi$ and $\chi' : T \longrightarrow o^{\times}$.}

Proof: We compute
$$
\matrix{
\Hom_{K[[B]]}(N_{\chi'},N_{\chi}^-) =
     \Hom_{K[[P]]}(K^{(\chi')},N_{\chi}^-)\hfill\cr\cr
\subseteq \Hom_{K[[U]]}(K,N_{\chi}^-) =
     \Hom_{K[[U]]}(K,K[[U]]) = 0\ .}
$$
The other vanishing follows by a completely symmetric computation.

\medskip

{\bf Theorem 4.4:}

{\it If $c(\chi) \not\in \Ndss_{0}$ then $M_{\chi}$ is a simple
$K[[G]]$-module.}

Proof: By our above results the decomposition $M_{\chi} \cong
N_{\chi} \oplus N_{w\chi}^-$ is a $K[[B]]$-invariant decomposition
into two nonisomorphic simple $K[[B]]$-modules. But as noted
already at the beginning it is not $K[[G]]$-invariant. Hence
$M_{\chi}$ must be a simple $K[[G]]$-module.

\medskip

The simple $K[[G]]$-modules which we have exhibited above are all
nonisomorphic as the following result implies.

\medskip

{\bf Proposition 4.5:}

{\it We have $\Hom_{K[[G]]} (M_{\chi'},M_{\chi}) = 0$ for any two
continuous characters $\chi \neq \chi' : T
\longrightarrow o^{\times}$.}

Proof: Because of Lemma 4.3 it is sufficient to show that
$\Hom_{K[[B]]} (N_{\chi'},N_{\chi}) = \Hom_{K[[B]]}
(N_{\chi'}^-,N_{\chi}^-) = 0$. Since the arguments are completely
symmetric we only discuss the vanishing of the first space. Making
as usual the identification $(\ast)$ we have
$$
\matrix{
\Hom_{K[[B]]}(N_{\chi'},N_{\chi}) =
     \Hom_{K[[P]]}(K^{(\chi')},N_{\chi})\hfill\cr\cr
= \{F \in K[[U^-]] : g(F) = \chi'(g)\cdot F\ \hbox{for any}\ g \in
P\}\ .}
$$
Assume now that there is a nonzero $F$ in this latter space. Since
any central matrix $g = \pmatrix{b & 0\cr 0 & b}$ in $T$ acts by
multiplication with $\chi(g)$ on $N_{\chi}$ it follows immediately
that $\chi$ and $\chi'$ have to coincide on those matrices. On the
other hand the action of an element $t_a \in T$ as described in
the proof of Prop. 4.1 gives rise to the equation
$$
\chi(t_a)\cdot F((1+x)^a - 1) = \chi'(t_a)\cdot F(x)
$$
between bounded power series over $K$. It was shown in the proof
of [ST1] Prop. 5.5 that this implies $c(\chi') - c(\chi) \in
2\Ndss_0$ and $F(x) = [{\rm log}(1+x)]^{(c(\chi')-c(\chi))/2}$.
Since the power series ${\rm log}(1+x)$ is not bounded we in fact
obtain $c(\chi') = c(\chi)$ and $F(x) = 1$. Going back to the
above equation it follows that $\chi(t_a) = \chi'(t_a)$. Hence we
have shown that the existence of a nonzero $F$ forces the
characters $\chi$ and $\chi'$ to coincide.

\medskip

To finish we briefly explain the dual picture. In the Banach space
$C(G,K)$ of all $K$-valued continuous functions on $G$ we have the
closed subspace
$$
\Ind_P^G(\chi):=\{ f\in C(G,K): f(gp)=
\chi(p^{-1})f(g) \ \hbox{for any}\ g\in G,
 p\in P\}\ .
$$
Via the left translation action this is a $K$-Banach space
representation of $G$ (a ''principal series'' representation). By
the interpretation of $K[[G]]$ as the space of bounded $K$-valued
measures on $G$ we have that $K[[G]]$ is the continuous dual of
$C(G,K)$. It easily follows that
$$
\Ind_P^G(\chi)' = M_{\chi^{-1}}\ .
$$
In particular, by Lemma 3.4, $\Ind_P^G(\chi)$ is an admissible
$G$-representation. As a consequence of Cor. 3.6 and Thm. 4.4 we
see that $\Ind_P^G(\chi)$ is topologically irreducible if $c(\chi)
\not\in -\Ndss_{0}$. This latter fact (for a slightly restricted
class of $\chi$) was proved in a direct and completely different
way in [Tru].

\bigskip

{\bf References}

\parindent=23truept

\ref{[AU]} Arnautov, V.I., Ursul, M.I.: On the uniqueness of topologies
for some constructions of rings and modules. Siberian Math. J. 36,
631-644 (1995)

\ref{[BH]} Ballister, P.N., Howson, S.: Note on Nakayama's Lemma
for Compact $\Lambda$-modules. Asian J. Math. 1, 224-229 (1997)

\ref{[B-CA]} Bourbaki, N.: Commutative Algebra. Paris: Hermann
1972

\ref{[B-GT]} Bourbaki, N.: General Topology. Berlin-Heidelberg-New
York: Sprin-ger 1989

\ref{[B-GAL]} Bourbaki, N.: Groupes et alg\`ebres de Lie, Chap. 2
et 3. Paris: Hermann 1972

\ref{[B-TVS]} Bourbaki, N.: Topological Vector Spaces.
Berlin-Heidelberg-New York: Springer 1987

\ref{[CR]} Curtis, C., Reiner, I.: Representation theory of finite
groups and associative algebras. New York-London: Wiley 1962

\ref{[SGA3]} Demazure, M., Grothendieck, A.: Sch\'emas en Groupes I.
Lect. Notes Math. 151. Berlin-Heidelberg-New York: Springer 1970

\ref{[Dia]} Diarra, B.: Sur quelques repr\'esentations
$p$-adiques de $\dZ$. Indagationes Math. 41, 481-493 (1979)

\ref{[Laz]} Lazard, M.: Groupes analytiques $\pfr$-adiques. Publ.
Math. IHES 26 (1965)

\ref{[Mon]} Monna, A.F.: Analyse non-archimedienne. Ergebnisse der
Math. 56. Berlin-Heidelberg-New York: Springer 1970

\ref{[Sch]} Schikhof, W.H.: A perfect duality between $p$-adic Banach
spaces and compactoids. Indag. Math. 6, 325-339 (1995)

\ref{[NFA]} Schneider, P.: Nichtarchimedische Funktionalanalysis.
Course at M\"un-ster (1997)

\ref{[ST1]} Schneider, P., Teitelbaum, J.: Locally analytic
distributions and $p$-adic representation theory, with
applications to $GL_2$. Preprintreihe des SFB 478, Heft 86,
M\"unster 1999

\ref{[ST2]} Schneider, P., Teitelbaum, J.: $U(\gfr)$-finite locally
analytic representations. Preprint 2000

\ref{[Ser]} Serre, J.-P.: Local Fields. Berlin-Heidelberg-New York: Springer
1979

\ref{[Tie]} van Tiel, J.: Espaces localement $K$-convexes I-III.
Indagationes Math. 27, 249-258, 259-272, 273-289 (1965)

\ref{[Tru]} Trusov, A.V.: Representations of the groups
$GL(2,\dZ)$ and $GL(2,\dQ)$ in spaces over non-archimedean fields.
Moscow Univ. Math. Bull. 36, 65-69 (1981)

\ref{[Was]} Washington, L.C.: Introduction to Cyclotomic Fields.
Berlin-Heidel-berg-New York: Springer 1982

\bigskip

\parindent=0pt

Peter Schneider\hfill\break Mathematisches Institut\hfill\break
Westf\"alische Wilhelms-Universit\"at M\"unster\hfill\break
Einsteinstr. 62\hfill\break D-48149 M\"unster, Germany\hfill\break
pschnei@math.uni-muenster.de\hfill\break
http://www.uni-muenster.de/math/u/schneider\hfill\break

\noindent
Jeremy Teitelbaum\hfill\break Department of Mathematics,
Statistics, and Computer Science (M/C 249)\hfill\break University
of Illinois at Chicago\hfill\break 851 S. Morgan St.\hfill\break
Chicago, IL 60607, USA\hfill\break jeremy@uic.edu\hfill\break
http://raphael.math.uic.edu/$\sim$jeremy\hfill\break

\end